%%%%%%%%%%%%%%%%%%%%%%%%%%%%%%%%%%%%%%%%%%%%%%%%%%%%%%%%%%%%%%%%%%%%%%%%%
% AMS-TeX file: Kaneko-Koike: Quasimodular forms as solutions to 
%a differential equation of hypergeometric type
% Version  2002.5.24 (written by Kaneko at OSU)
% revised  2002.6.4 (after correction by Koike)
% revised  2002.7.26 (after referee's report)
%%%%%%%%%%%%%%%%%%%%%%%%%%%%%%%%%%%%%%%%%%%%%%%%%%%%%%%%%%%%%%%%%%%%%%%%%

\input amstex.tex
\documentstyle{amsppt}
\mag=\magstep1
\pagewidth{160truemm}
\pageheight{230truemm}%A4size
\topmatter
\title \nofrills \bf {Quasimodular solutions of a differential
equation of hypergeometric type} \endtitle
\rightheadtext{Quasimodular forms}
\leftheadtext{Masanobu Kaneko and Masao Koike}
\author {Masanobu Kaneko and Masao Koike} \endauthor
\address {Graduate School of Mathematics, 
Kyushu University 33, Fukuoka 812-8581, Japan}\endaddress
\email {mkaneko\@math.kyushu-u.ac.jp} \endemail
\address {Graduate School of Mathematics, Kyushu University,
Ropponmatu, Fukuoka 810-8560, Japan} 
\endaddress 
\email{koike\@math.kyushu-u.ac.jp}\endemail
\endtopmatter 
\document
\tolerance=500 
\nologo
\baselineskip=1.1\baselineskip
%\pageno=1
%%%%%%%%%%%%%%%%%%%% 

\def\section#1{\heading{#1}\endheading}

 \def\={\;=\;} \def\:{\;:=\;}

  \def\C{\Bbb C}  \def\D{\Delta}
\def\d{\delta} \def\d2{\D_{2A}}  \def\eq{$(\#)_k$ } \def\eqd{$(\#')_k$ }
    \def\G{\Gamma}
\def\g{\text{SL}_2(\Z)}  \def\h{\frak{H}}  \def\i{^{-1}}
\def\inf{\infty}   \def\j1{j(\tau)} \def\j2{j_{(2)}(\tau)}
\def\j3{j_{(3)}(\tau)} \def\j4{j_{(4)}(\tau)} 
\def\l{\lambda}    
     \def\Q{\Bbb Q} 
    \def\t{\tau} \def\th{\vartheta}

 \def\Z{\Bbb Z}

\section{\S1. Introduction and Main Theorem}

In our previous paper [2], we studied further the solutions of the
following differential equation in the upper half-plane $\h$ which was
originally found and studied in [4] in connection with the arithmetic
of supersingular elliptic curves;
$$f''(\tau)-\frac{k+1}6 E_2(\tau)f'(\tau)+\frac{k(k+1)}{12}
E'_2(\tau)f(\tau)=0.$$
Here, $k$ is an integer or half an integer, the
symbol ${}'$ denotes  the differentiation $(2\pi i)\i d/d\tau=q\cdot
d/dq\ (q=e^{2\pi i \tau})$, and $E_2(\tau)$ is the ``quasimodular''
Eisenstein series of weight $2$ for the full modular group $\g$:
$$E_2(\tau)=1-24\sum_{n=1}^\inf\bigl(\sum_{d|n} d\bigr)q^n.$$

Let $p\ge5$ be a prime number and $F_{p-1}(\tau)$ be the solution of
the above differential equation for $k=p-1$ which is modular on $\g$
(such a solution exists and is unique up to a scalar multiple). For
any zero $\tau_0$ in $\h$ of the form $F_{p-1}(\t)$, the value of the
$j$-function at $\t_0$ is algebraic and its reduction modulo (an
extension of) $p$ is a supersingular $j$-invariant of characteristic
$p$, and conversely, all the supersingular $j$-invariants are obtained
in this way from the single solution $F_{p-1}(\t)$ with suitable
choices of $\t_0$.  This is the arithmetic connection that motivated our
study of the differential equation.  

Various modular forms on $\g$ and its subgroups were obtained in [2]
as solutions to this differential equation, the groups depending on
the choice of $k$. Every modular solution is expressed in terms of 
a hypergeometric polynomial in a suitable modular function (hence the
``hypergeometric type'' in the title of the paper), also 
depending on the choice of $k$. For instance, if $k\equiv0, 4\bmod 12$,
we have a modular solution 
$$E_4(\tau)^{\frac{k}4}
F(-\frac{k}{12},-\frac{k-4}{12},-\frac{k-5}6;\frac{1728}{j(\t)}),$$
where $$F(a,b,c;x)=\sum_{n=0}^\inf \frac{(a)_n(b)_n}
{(c)_n}\frac{x^n}{n!},\qquad (a)_n=a(a+1)\cdots (a+n-1)$$
is the Gauss
hypergeometric series (which becomes a polynomial when $a$ or $b$ is a
negative integer, which is the case here), $E_4(\t)$ the Eisenstein
series of weight $4$ on $\g$, and $j(\t)$ the elliptic modular
invariant.

In addition to the modular solutions, quite remarkable  was an
occurrence of a {\it quasimodular form}, not of weight $k$ as in the
modular case but of weight $k+1$. In the present paper, we give
another supply of examples of quasimodular forms as solutions  to an
analogous differential equation attached to the group $\Gamma_0^*(2)$,
which is {\it not} contained in $\g$;  
$$\G_0^*(2)=\left\langle \G_0(2),\ \pmatrix 0 & -1 \\ 2 & 0 \endpmatrix
\right\rangle$$ where $$\G_0(2)=\left\{\pmatrix a&b\\c&d
    \endpmatrix\in\g  |\; c\equiv 0\mod 2\right\}.$$
($\Gamma_0^*(2)$ is the triangular  group ``$2A$'' in the notation of
Conway-Norton [1].)

Let $$E_{2A}(\tau):=(E_2(\tau)+2E_2(2\tau))/3=1 - 8q - 40q^2 - 32q^3
-\cdots$$
be the quasimodular form of weight $2$ on $\G_0^*(2)$ which
is the logarithmic derivative of the form
$$\Delta_{2A}(\tau):=\eta(\tau)^8\eta(2\tau)^8=q - 8 q^2 + 12 q^3  +
64 q^4 - \cdots$$
of  weight $8$ on $\G_0^*(2)$;
$E_{2A}(\tau)=\Delta_{2A}'(\tau)/\Delta_{2A}(\tau)$,   an analogous
situation in the $SL_2(\Z)$ case where $E_2(\tau)$ is the logarithmic
derivative of the Ramanujan $\Delta(\tau)$.   Consider the following
differential equation;
$$(\#)_k\qquad f''(\tau)-\frac{k+1}4
E_{2A}(\tau)f'(\tau)+\frac{k(k+1)}{8} E'_{2A}(\tau)f(\tau)=0.$$
Solutions which are modular on the group $\G_0^*(2)$ and its subgroups
were studied in [6, 7]. In particular, when $k$ is a non-negative
integer  congruent to $0$ or $6$ modulo $8$, the equation \eq has a
one dimensional space of solutions which are modular on the group
$\G_0^*(2)$ itself.  We note here that the equation \eq has a
characterization by the invariance of the space of solutions under the
action of $\G_0^*(2)$, similar to  the previous case for $\g$, owing to
the fact that there is no holomorphic modular form of weight $2$ on
$\G_0^*(2)$ (see [5] and [2, \S5]).  By a general theory of ordinary
differential equations, we see that the equation \eq has a
quasimodular solution (which, since its transformation under
$\tau\rightarrow -1/2\tau$ is also a solution,  inevitably gives a
solution having $\log q$ term in the expansion at $q=0$) only when $k$
is a positive integer congruent to $3$ modulo $4$.

In the following, we show there indeed exists a quasimodular solution
in this case and describe explicitly the solution in terms of a
certain orthogonal polynomials.  First we need to develop some
notations.  Put
$$\align  C(\t) &:= 2E_2(2\t)-E_2(\t)\\&=1+24\sum_{n=1}^\inf \bigl(
\sum_{d|n\atop d:\text{odd}} d\; \bigr)q^n= 1+24q+24q^2+96q^3+\cdots,\\
D(\t) &:= \frac{\eta(2\t)^{16}}{\eta(\t)^8}= \sum_{n=1}^\inf
\bigl(\sum_{d|n\atop d:\text{odd}} (n/d)^3\bigr)
q^n=q+8q^2+28q^3+64q^4+\cdots, \endalign $$
where
$$\eta(\t)=q^{\frac1{24}}\prod_{n=1}^\inf (1-q^n) =q^\frac{1}{24} -
q^\frac{25}{24}- q^\frac{49}{24} + q^\frac{121}{24}+ \cdots$$
is the
Dedekind eta function.  The functions $C(\t)$ and $D(\t)$ are modular
forms of respective weights $2$ and $4$ on the group
$\G_0(2)$ (=``$2B$'') and the graded ring of modular  forms of integral
weights on $\G_0(2)$ is generated by these $C(\t)$ and $D(\t)$.
Recall that (see [3])  an element of degree $k$ in the graded ring
$\C[E_2(\t),C(\t),D(\t)]$, where the generators $E_2(\t), C(\t),
D(\t)$ have degrees $2,\, 2$, and $4$ respectively, is referred to as
a quasimodular form of weight $k$ (on $\Gamma_0(2)$). Incidentally,
the graded ring of modular  forms of integral weights on $\G_0^*(2)$ is
generated by three elements $C(\t)^2=(E_4(\t)+4E_4(2\t))/5$,
$C(\t)^3-128C(\t)D(\t)=(E_6(\t)+8E_6(2\t))/9$, and $\D_{2A}(\t)$ of
respective weights $4\,,6\,,8$, of  which $C(\t)^2$ and $\D_{2A}(\t)$
generate freely the subring consisting forms of weight being multiple
of $4$, and the whole  space as a graded module is generated over this
ring  by  $C(\t)^3-128C(\t)D(\t)$.
 
Now define a sequence of polynomials $P_n(x)\ (n=0,1,2,\dots)$ by
$$P_0(x)=1,\  P_1(x)=x,\  \ P_{n+1}(x)=xP_n(x)+\l_nP_{n-1}(x)\quad
(n=1,2,\dots)$$
where $$\l_n=4\frac{(4n+1)(4n+3)}{n(n+1)}.$$
First few
examples are $$P_2(x)=x^2+70,\  P_3(x)=x^3+136x,\
P_4(x)=x^4+201x^2+4550,\dots.$$
The $P_n(x)$ is  even or odd
polynomial according as $n$ is even or odd. We also define a second
series of polynomials $Q_n(x)$ by the same recursion  (with different
initial values):
$$Q_0(x)=0,\  Q_1(x)=1,\ \ Q_{n+1}(x)=xQ_n(x)+\l_nQ_{n-1}(x)\quad
(n=1,2,\dots),$$
a couple of examples being $$Q_2(x)=x,\  Q_3(x)=x^2+66,\ 
Q_4(x)= x^3+131x,\dots.$$
The $Q_n(x)$ has opposite parity: It is even
if $n$  is odd and odd if $n$ is even.

Put $G(\t)=C(\t)^2-128D(\t)\ (=(4E_4(2\t)-E_4(\t))/3)$.

\proclaim{Theorem} Let $k=4n+3\ \ (n=0,1,2,\dots)$. The following 
quasimodular form of weight $k+1$ on $\Gamma_0(2)$ is a solution of $(\#)_k:$
$$\sqrt{\D_{2A}(\t)}^nP_n\bigl(\frac{G(\t)}{\sqrt{\D_{2A}(\t)}}\bigr)
\frac{C'(\t)}{24}
-\sqrt{\D_{2A}(\t)}^{n+1}Q_n\bigl(\frac{G(\t)}{\sqrt{\D_{2A}(\t)}}\bigr).$$
\endproclaim

\example{Remark} The appearance of the square root
$\sqrt{\D_{2A}(\t)}$ in the formula is superficial  because of the
parities of $P_n(x)$ and $Q_n(x)$, that is, the form is actually an
element in $\Q[ E_2(\t),C(\t),D(\t)]$, by noting
$\D_{2A}(\t)=D(\t)(C(\t)^2- 64D(\t))$ and
$C'(\t)=(E_2(\t)C(\t)-C(\t)^2)/6+ 32D(\t)$. The form does not belong
to $\G_0^*(2)$.  \endexample
 
\section{\S2. Proof of Theorem}

It is convenient to  introduce the operator $\th_k$ defined by
$$\th_k(f)(\t)=f'(\t)-\frac{k}{8}E_{2A}(\t)f(\t).$$
By the quasimodular property of $E_2(\t)$ or the fact that $E_{2A}(\t)$ is
the logarithmic derivative of $\D_{2A}(\t)$, we have the
transformation formulas
$$E_{2A}\left(\frac{a\t+b}{c\t+d}\right)=(c\t+d)^2E_{2A}(\t)+\frac4{\pi i}
c(c\t+d) \qquad (\pmatrix a & b \\ c & d \endpmatrix \in\G_0(2))$$ and
$$E_{2A}\left(-\frac1{2\t}\right)=2\t^2E_{2A}(\t)+\frac8{\pi i}\t.$$
From these we see that if $f$ is modular of weight $k$ on a subgroup
of $\G_0^*(2)$,  then $\th_k(f)$  is modular of weight $k+2$ on the same
group.  If $f$ and $g$ have weights $k$ and $l$, the Leibniz rule
$$\th_{k+l}(fg)=\th_k(f)g+f\th_l(g)$$
holds. We sometimes drop the
suffix  of the operator $\th_k$ when the weights of modular forms we
consider are clear. With this operator, the equation \eq can be
rewritten as
$$(\#')_k\qquad\th_{k+2}\th_k(f)(\t)=\frac{k(k+2)}{64}C(\t)^2f(\t),
\hskip60pt$$ (use $E_{2A}'(\t)=(E_{2A}(\t)^2-C(\t)^2)/8$).

Denote the form in the theorem by $F_k(\t)$.  We first establish the
recurrence relation (note $n=(k-3)/4$):
$$ F_{k+4}(\t)=G(\t) F_k(\t)+\l_{n}\D_{2A}(\t)F_{k-4}(\t). \tag1$$
This is a consequence of the recursion of $P_n$ and $Q_n$, namely
(we often omit the variable $\t$ hereafter)
$$
\align &GF_{k}+\l_{n}\D_{2A} F_{k-4}\\
&=G\biggl(\sqrt{\D_{2A}}^{n}P_{n}\bigl(\frac{G}{\sqrt{\D_{2A}}}\bigr)
\frac{C'}{24}
-\sqrt{q\D_{2A}}^{n+1}Q_{n}\bigl(\frac{G}{\sqrt{\D_{2A}}}\bigr)\biggr)\\
&+\l_{n}\D_{2A}\biggl(
\sqrt{\D_{2A}}^{n-1}P_{n-1}\bigl(\frac{G}{\sqrt{\D_{2A}}}\bigr)\frac{C'}{24}
-\sqrt{\D_{2A}}^{n}Q_{n-1}\bigl(\frac{G}{\sqrt{\D_{2A}}}\bigr)\biggr)\\
&=\sqrt{\D_{2A}}^{n+1}\biggl(\frac{G}{\sqrt{\D_{2A}}}
P_{n}\bigl(\frac{G}{\sqrt{\D_{2A}}}\bigr)+\l_{n}
P_{n-1}\bigl(\frac{G}{\sqrt{\D_{2A}}}\bigr)\biggr)\frac{C'}{24}\\
&-\sqrt{\D_{2A}}^{n+2}\biggl(\frac{G}{\sqrt{\D_{2A}}}
Q_{n}\bigl(\frac{G}{\sqrt{\D_{2A}}}\bigr)+\l_{n}
Q_{n-1}\bigl(\frac{G}{\sqrt{\D_{2A}}}\bigr)\biggr)\\
&=\sqrt{\D_{2A}}^{n+1}P_{n+1}\bigl(\frac{G}{\sqrt{\D_{2A}}}\bigr)
\frac{C'}{24}
-\sqrt{\D_{2A}}^{n+2}Q_{n+1}\bigl(\frac{G}{\sqrt{\D_{2A}}}\bigr)\\
&=F_{k+4}. \endalign $$

Now we prove by induction that the $F_k(\t)$ satisfies the equation
$(\#')_k$.  We can check the cases $k=3$ and $7$ directly. Assume
$F_{k-4}$ and  $F_k$ satisfy $(\#')_{k-4}$ and \eqd respectively.
Then by using (1) and the formulas
$$\th(C)=-\frac14G,\quad \th(G)=-\frac12C^3,\quad \th(\D_{2A})=0$$
we have
$$\align \th^2(F_k)&=\th\bigl(\th(F_k)G-\frac12C^3F_k\bigr)+\l_n\D_{2A}
\th^2(F_{k-4})\\
&=\th^2(F_k)G-\frac12\th(F_k)C^3+\frac38C^2GF_k-\frac12C^3\th(F_k)
+\l_n\d2\th^2(F_{k-4})\\
&=\frac{k(k+2)}{64}C^2GF_k-C^3\th(F_k)+\frac38C^2GF_k+\frac{(k-4)(k-2)}{64}
\l_n\d2 C^2F_{k-4}\\
&= \frac{k^2+2k+24}{64}C^2GF_k+\frac{(k-4)(k-2)}{64}\l_n\d2 C^2F_{k-4}
-C^3\th(F_k). \endalign $$
Hence we find 
$$\align &\th^2(F_{k+4})-\frac{(k+4)(k+6)}{64}C^2F_{k+4}\\
&=\biggl(\frac{k^2+2k+24}{64}-\frac{(k+4)(k+6)}{64}\biggr)C^2GF_k\\
&+
\biggl(\frac{(k-4)(k-2)}{64}-\frac{(k+4)(k+6)}{64}\biggr)\l_n\d2 C^2F_{k-4}\\
&=-C^2\biggl(\frac{k}8GF_k+C\th(F_k)+\frac{k+1}4\l_n\d2 F_{k-4}\bigr).
\endalign $$
The proof of the theorem therefore boils down to show the
equation
$$\frac{k}8GF_k+C\th(F_k)=-\frac{k+1}4\l_n\d2 F_{k-4}.$$
For this we also proceed by induction. For $k=7$ the equation is  checked
directly. Assuming that this is valid for $k$, we have
$$F_{k+4}=GF_k+\l_n\d2 F_{k-4}=\frac1{2(k+1)}\bigl((k+2)GF_k-8C\th(F_k)\bigr)$$
and
$$\align &\frac{k+4}8GF_{k+4}+C\th(F_{k+4})\\
&=\frac{k+4}{16(k+1)}G\bigl((k+2)GF_k-8C\th(F_k)\bigr)\\&+
\frac1{2(k+1)}C\bigl(-\frac12(k+2)C^3F_k+(k+2)G\th(F_k)+2G\th(F_k)-8C\th^2(F_k)
\bigr)\\
&=\frac{(k+2)(k+4)}{16(k+1)}(G^2-C^4)F_k\\
&=-\frac{k+5}4\l_{n+1}\d2 F_k. \endalign $$
Here we have used the
(previous) induction assumption that $F_k$  satisfies \eqd  and the
relation $G^2-C^4=-256\d2$. This completes our proof.

\enddocument